# Approximation of the derivatives of the logarithm of the Riemann zeta-function in the critical strip

Sergey K. Sekatskii[a] and Stefano Beltraminelli[b]

[a]*Laboratoire de Physique de la Matière Vivante, IPHYS, Ecole Polytechnique Fédérale de Lausanne, BSP, CH 1015 Lausanne, Switzerland* E-mail : Serguei.Sekatski@epfl.ch
[b]*CERFIM, Research Center for Mathematics and Physics, PO Box 1132, 6600 Locarno, Switzerland*

Recently, we have established the generalized Li's criterion equivalent to the Riemann hypothesis, viz. demonstrated that the sums over all non-trivial Riemann function zeroes $k_{n,a} = \Sigma_\rho (1 - \left(\frac{\rho - a}{\rho + a - 1}\right)^n)$ for any real $a$ not equal to ½ are non-negative if and only if the Riemann hypothesis holds true, and proved the relation $k_{n,a} = \frac{n(1-2a)}{(n-1)!}\frac{d^n}{dz^n}((z-a)^{n-1}\ln(\xi(z)))|_{z=1-a}$. Assuming that the function $\varsigma(s)$ is non-vanishing for $\operatorname{Re} s > 1/2 + \Delta$, where real $0 < \Delta < 1/2$, using this relation together with the functional equation for the $\xi$-function and the explicit formula of Weil, we prove that in these conditions for $n=1, 2, 3…$ and an arbitrary complex $a$ with $1 > \operatorname{Re} a \geq 1/2 + \Delta + \delta_0$, where $\delta_0$ is an arbitrary small fixed positive number, one has

$$\frac{d^n}{ds^n}\ln \varsigma(s)|_{s=a} = \sum_{m \leq N} \frac{(-1)^{n-1}\Lambda(m)\ln^{n-1} m}{m^a} + \int_0^N x^{-a}\ln^{n-1} x\, dx + O(N^{1/2+\Delta-a}\ln^{n-1} N).$$

In particular, $\frac{\varsigma'}{\varsigma}(a) = -\sum_{m \leq N}\frac{\Lambda(m)}{m^a} + \frac{N^{1-a}}{1-a} + O(N^{1/2+\Delta-a})$. Numerical verifications of these equalities are also presented.

**Introduction**

In a recent paper [1], we have established the generalized Li's criterion equivalent to the Riemann hypothesis and first discovered in [2] (see e.g. [3] for general discussion of properties of the Riemann zeta-function), as well as the closely related generalized Bombieri – Lagarias theorem concerning the location of zeroes of certain complex number multisets [4]. Namely, we have demonstrated that the sums $k_{n,a} = \Sigma_\rho (1 - \left(\frac{\rho - a}{\rho + a - 1}\right)^n)$, $n$=1, 2, 3..., over the non-trivial Riemann zeta-function zeroes $\rho = \sigma + iT$ taking into account their multiplicities, for any real $a$ not equal to ½ are non-negative if and only if the Riemann hypothesis holds true. (Throughout the paper we understand the sums over non-trivial Riemann zeroes as $\sum_\rho f(\rho) = \lim_{X \to \infty,\ |T|<X} \sum_\rho f(\rho)$ and will not mention this any more). We also established the relation between these sums and certain derivatives of the Riemann $\varsigma$-function: $k_{n,a} = \frac{n(1-2a)}{(n-1)!} \frac{d^n}{dz^n}((z-a)^{n-1} \ln(\xi(z)))|_{z=1-a}$ and gave an arithmetic interpretation of these sums.

These results, as well as the arXiv submissions [5, 6] of one of the authors (S.K.S), are used in this Note, which aim is to establish the following theorem 1 and then, which seems more interesting for us, its corrolary.

**Theorem 1.** *Assume the Riemann function $\varsigma(s)$ is non-vanishing for $\operatorname{Re} s > 1/2 + \Delta$ where real $0 < \Delta < 1/2$. Then for n=1, 2, 3... and an arbitrary complex a with $1 > \operatorname{Re} a \geq 1/2 + \Delta + \delta_0$, where $\delta_0$ is an arbitrary small fixed positive number, we have:*

$$\sum_{\rho}(1-\left(\frac{\rho-a}{\rho+a-1}\right)^n) = \sum_{\rho}(1-\left(\frac{\rho+a-1}{\rho-a}\right)^n) = 2-(-1+\frac{1}{a})^n-(-1+\frac{1}{1-a})^n$$

$$\sum_{j=1}^{n} C_n^j (2a-1)^j \left\{ \frac{(-1)^j}{(j-1)!} \lim_{N\to\infty} (\sum_{m\leq N} \frac{\Lambda(m)\ln^{j-1} m}{m^a} - \int_0^N x^{-a}\ln^{j-1} x\, dx) \right\} + \qquad (1)$$

$$+\frac{n}{2}(2a-1)(\psi(a/2)-\ln\pi) + \sum_{j=2}^{n} C_n^j (-1)^j 2^{-j}(2a-1)^j \varsigma(j, a/2)$$

Here $\Lambda(n)$ is van Mangoldt function and $\varsigma(s, a) := \sum_{m=0}^{\infty} \frac{1}{(m+a)^s}$ is Hurwitz zeta-function.

This immediately implies the following corrolary.

**Corrolary 1**. *In conditions of Theorem 1, we have for any integer n=1, 2, 3...:*

$$\frac{d^n}{ds^n}\ln\varsigma(s)|_{s=a} = \sum_{m\leq N} \frac{(-1)^{n-1}\Lambda(m)\ln^{n-1} m}{m^a} + \int_0^N x^{-a}\ln^{n-1} x\, dx + O(N^{1/2+\Delta-a}\ln^{n-1} N) \qquad (2).$$

*In particular,*

$$\frac{\varsigma'}{\varsigma}(a) = -\sum_{m\leq N}\frac{\Lambda(m)}{m^a} + \frac{N^{1-a}}{1-a} + O(N^{1/2+\Delta-a}) \qquad (3).$$

*Assuming the validity of the Riemann hypothesis, we evidently have*

$$\frac{d^n}{ds^n}\ln\varsigma(s)|_{s=a} = \sum_{m\leq N}\frac{(-1)^{n-1}\Lambda(m)\ln^{n-1} m}{m^a} + \int_0^N x^{-a}\ln^{n-1} x\, dx + O(N^{1/2-a}\ln^{n-1} N).$$

Driven by an analogy with the quite known $\varsigma(a) = \sum_{m\leq N}\frac{1}{m^a} - \frac{N^{1-a}}{1-a} + O(N^{-\text{Re}\,a})$, unconditionally valid in the critical strip, Theorem 4.11 of [3] (additional conditions on *Im(a)* must be added here, see [3]), we will name these relations "conditional approximation of the derivatives of the logarithm of Riemann zeta-function in the critical strip".

**Remark 1.** This is well known that for Re$s>1$, $\dfrac{\varsigma'(s)}{\varsigma(s)} = -\sum_{m=1}^{\infty} \dfrac{\Lambda(m)}{m^s}$ [3].

**Remark 2**. Certainly, in eqs. (1, 2) and other places in the text we can also use the well-known elementary integral

$$\int_0^N x^{-a} \ln^{j-1} x\, dx = (-1)^{j-1} \frac{d^{j-1}}{da^{j-1}} \int_0^N x^{-a} dx = (-1)^{j-1} \frac{d^{j-1}}{da^{j-1}} \frac{N^{1-a}}{1-a} =$$

$$\sum_{l=1}^{j} C_{j-1}^{l-1} \frac{(-1)^{l+j} N^{1-a} (l-1)!}{(1-a)^l} \ln^{j-l} N = \sum_{l=1}^{j} \frac{(-1)^{l+j} (j-1)! N^{1-a}}{(1-a)^l (j-l)!} \ln^{j-l} N \qquad (4).$$

We do not do this here to simplify the notation.

## 1. An arithmetic interpretation

For completeness, we start with the repetition of the material given in [1].

For suitable function $f$, Mellin transform is defined as $\hat{f}(s) = \int_0^{\infty} f(x) x^{s-1} dx$ while inverse Mellin transform formula is $f(x) = \dfrac{1}{2\pi i} \int_{\text{Re}\, s=c} \hat{f}(s) x^{-s} ds$ with an appropriate value of $c$.

**Lemma 1.** *For $n=1, 2, 3...$, and an arbitrary complex number $a$ the inverse Mellin transform of the function $k_{n,a}(s) = 1 - \left(1 - \dfrac{2a-1}{s+a-1}\right)^n$ is*

$$g_{n,a}(x) = P_{n,a}(x) \qquad \text{if } 0 < x < 1$$

$$g_{n,a}(x) = \frac{n}{2}(2a-1) \quad \text{if } x = 1 \qquad (5)$$

$$g_{n,a}(x) = 0 \qquad \text{if } x > 1$$

*where* $P_{n,a}(x) = x^{a-1} \sum_{j=1}^{n} C_n^j \dfrac{(2a-1)^j \ln^{j-1} x}{(j-1)!}$; $C_n^j = \dfrac{n!}{j!(n-j)!}$ *is a binomial coefficient.*

*Proof.* We have for Re$(s+a)>1$:

$$\sum_{j=1}^{n} C_n^j \frac{(2a-1)^j}{(j-1)!} \int_0^1 (\ln^{j-1} x) x^{s+a-2} dx = \sum_{j=1}^{n} C_n^j \frac{(2a-1)^j}{(j-1)!} \frac{d^{j-1}}{ds^{j-1}} \int_0^1 x^{s+a-2} dx =$$

$$\sum_{j=1}^{n} C_n^j \frac{(2a-1)^j (-1)^{j-1}}{(s+a-1)^j} = 1 - \left(1 - \frac{2a-1}{s+a-1}\right)^n$$

If $a$ is an arbitrary complex number with $\mathrm{Re}\, a > 1$, for the function $g_{n,a}(x)$ we can apply the so called Explicit Formula of Weil, see [4, 7, 8], which is, as given in [4]:

$$\sum_{\rho} \hat{f}(\rho) = \int_0^\infty f(x)dx + \int_0^\infty \tilde{f}(x)dx - \sum_{n=1}^{\infty} \Lambda(n)(f(n) + \tilde{f}(n))$$

$$- (\ln \pi + \gamma) \cdot f(1) - \int_1^\infty \left\{ f(x) + \tilde{f}(x) - \frac{2}{x^2} f(1) \right\} \frac{x dx}{x^2 - 1} \quad (6)$$

and we define $\tilde{f}(x) := \frac{1}{x} f(\frac{1}{x})$, thus in our case the function

$$\tilde{P}_{n,a}(x) = x^{-a} \sum_{j=1}^{n} C_n^j \frac{(-1)^{j-1}(2a-1)^j \ln^{j-1} x}{(j-1)!}$$ should be used whenever appropriate.

Certainly, $P_{n,a}(x) = x^{a-1}(2a-1) L_{n-1}^1(-(2a-1) \cdot \ln x)$ and

$\tilde{P}_{n,a}(x) = x^{-a}(2a-1) L_{n-1}^1((2a-1) \cdot \ln x)$, where $L_n(x) = \sum_{j=0}^{n} C_n^j \frac{(-1)^j x^j}{j!}$ and

$\frac{d}{dx} L_n(x) = \sum_{j=1}^{n} C_n^j \frac{(-1)^j x^{j-1}}{(j-1)!} = -L_{n-1}^1(x)$ are generalized Laguerre polynomials [9], cf. [10].

This is easy to check that the function $g_{n,a}(x)$ do possess the necessary properties in a sense of continuity and asymptotic (in particular, for some positive $\delta$, $g_{n,a}(x) = O(x^\delta)$ as $x \to 0+$) for eq. (6) to be true [4, 7, 8].

Such an application gives

$$\sum_{\rho}(1-\left(\frac{\rho-a}{\rho+a-1}\right))^n = \sum_{\rho}(1-\left(\frac{\rho+a-1}{\rho-a}\right))^n =$$

$$\sum_{j=1}^{n} C_n^j \frac{(2a-1)^j}{(j-1)!} \{\int_0^1 x^{a-1} \ln^{j-1} x dx + (-1)^{j-1}\int_1^{\infty} x^{-a} \ln^{j-1} x dx - (-1)^{j-1}\sum_{m=1}^{\infty}\frac{\Lambda(m)\ln^{j-1} m}{m^a}\}$$

$$-\frac{n}{2}(2a-1)(\ln\pi+\gamma) - \int_1^{\infty}\{\sum_{j=1}^{n} C_n^j \frac{(-1)^{j-1}(2a-1)^{j-1}}{(j-1)!} x^{-a}\ln^{j-1} x - \frac{n}{x^2}(2a-1)\}\frac{xdx}{x^2-1}$$ (7).

Now, in the second and third integrals in the r.h.s. of (7) we make a variable transform $x$ to $1/x$, after what these integrals take the forms

$$I_2 = \int_0^1 x^{a-2}\ln^{j-1}(x)dx \text{ and } I_3 = \int_0^1\{\sum_{j=2}^{n} C_n^j \frac{\ln^{j-1} x}{(j-1)!}(2a-1)^j x^{a-1} + n(2a-1)(x^{a-1}-x)\}\frac{dx}{1-x^2}.$$

The first two integrals are handled by virtue of an example 4.272.6 of GR book [9]: $\int_0^1 \ln^{\mu-1}(1/x)x^{\nu-1}dx = \frac{1}{\nu^{\mu}}\Gamma(\mu)$; $\text{Re}\mu > 0$ and $\text{Re}\nu > 0$. Adopting for our case, we get $\int_0^1 \ln^{j-1}(x)x^{a-1}dx = \frac{(-1)^{j-1}}{a^j}(j-1)!$, $\int_0^1 \ln^{j-1}(x)x^{a-2}dx = \frac{(-1)^{j-1}}{(a-1)^j}(j-1)!$.

The "second part" of the third integral $I_3$ is, by virtue of an example 3.244.3 of the book [11], equal to

$$I_{32} = n(2a-1)\int_0^1 \frac{x^{a-1}-x}{1-x^2}dx = -\frac{n}{2}(2a-1)(\gamma+\psi(a/2)); \text{ here } \gamma = 0.572... \text{ is}$$

Euler – Mascheroni constant and $\psi$ is a digamma function. In the first part of this integral we make the variable change $x=\exp(-t)$:

$$I_{31} = \int_0^1\sum_{j=2}^{n} C_n^j \frac{\ln^{j-1} x}{(j-1)!}(2a-1)^j x^{a-1}\frac{dx}{1-x^2} = \sum_{j=2}^{n} C_n^j (-1)^{j-1}\frac{(2a-1)^j}{(j-1)!}\int_0^{\infty} t^{j-1}\frac{e^{-at}}{1-e^{-2t}}dt.$$

Applying Taylor expansion $(1-e^{-2t})^{-1} = 1+e^{-2t}+e^{-4t}+e^{-6t}+...$ we get further

$$I_{31} = \sum_{j=2}^{n} C_n^j(-1)^{j-1}\frac{(2a-1)^j}{(j-1)!}\sum_{m=0}^{\infty}\frac{(j-1)!}{(2m+a)^j} = \sum_{j=2}^{n} C_n^j(-1)^{j-1}2^{-j}(2a-1)^j\varsigma(j,a/2).$$

Using the relations

$$\sum_{j=1}^{n} C_n^j (-1)^{j-1} (2a-1)^j a^{-j} = 1 - \sum_{j=0}^{n} C_n^j (-1)^j \left(\frac{2a-1}{a}\right)^j = 1 - (-1+\frac{1}{a})^n,$$

$$\sum_{j=1}^{n} C_n^j (-1)^{j-1} (2a-1)^j (a-1)^{-j} = 1 - (-1+\frac{1}{1-a})^n,$$ and collecting everything together

we have proven the following theorem [1]:

**Theorem 2** *For n=1, 2, 3... and an arbitrary complex a with* $\operatorname{Re} a > 1$
*we have*

$$\sum_{\rho} (1 - \left(\frac{\rho-a}{\rho+a-1}\right)^n) = \sum_{\rho} (1 - \left(\frac{\rho+a-1}{\rho-a}\right)^n) = 2 - (-1+\frac{1}{a})^n - (-1+\frac{1}{1-a})^n +$$

$$\sum_{j=1}^{n} C_n^j (2a-1)^j \frac{(-1)^j}{(j-1)!} \sum_{m=1}^{\infty} \frac{\Lambda(m) \ln^{j-1} m}{m^a} + \frac{n}{2}(2a-1)(\psi(a/2) - \ln \pi) + \qquad (8).$$

$$\sum_{j=2}^{n} C_n^j (-1)^j 2^{-j} (2a-1)^j \varsigma(j, \, a/2)$$

Now we proceed to the proof of the Theorem 1.

*Proof.* In the consitions of the theorem the functions $g_{n,a}(x)$ can not be used because they do not possess the necessary asymptotic for *x* tending to zero. We need to introduce "truncated" functions $g_{n,a,\varepsilon}(x)$, $0 < \varepsilon < 1$, defined as follows:

$$g_{n,a,\varepsilon}(x) = g_{n,a}(x) \quad \text{if } \varepsilon < x \leq \infty$$

$$g_{n,a,\varepsilon}(x) = \frac{1}{2} g_{n,a}(\varepsilon) \text{ if } x = \varepsilon$$

$$g_{n,a,\varepsilon}(x) = 0 \qquad \text{if } x < \varepsilon.$$

The repetition of the calculations presented above during the proof of Theorem 2 with these truncated functions furnishes eq. (1): the terms containing the limits comes from the term $\int_0^{\infty} \tilde{f}(x) dx$ in eq. (5), which takes the form $\int_{\varepsilon}^{1} x^{a-2} \ln^{j-1} x \, dx$, and from the term $\sum_{n=1}^{\infty} \Lambda(n) \tilde{f}(n)$. (Initially, we obtain the limit $\lim_{N\to\infty}(\sum_{m\leq N}\frac{(-1)^{n-1}\Lambda(m)\ln^{n-1} m}{m^a}-\int_1^N x^{-a}\ln^{n-1} xdx)$, which we transform to

$\lim_{N\to\infty}(\sum_{m\leq N}\frac{(-1)^{n-1}\Lambda(m)\ln^{n-1} m}{m^a}-\int_0^N x^{-a}\ln^{n-1} xdx)$ using $\int_0^1 x^{-a}\ln^{j-1} xdx = \frac{(-1)^{j-1}(j-1)!}{(1-a)^{j-1}}$,

sf. the proof of Theorem 2).

Thus it rests to show that $\lim_{\varepsilon\to 0+}\sum_\rho \hat{g}_{n,a,\varepsilon}(\rho)=\sum_\rho \hat{g}_{n,a}(\rho)$. For small enough $\varepsilon$ we have $|\hat{g}_n(s)-\hat{g}_{n,\varepsilon}(s)|=|\int_0^\varepsilon T_n(\ln x)x^{a+s-2}dx|<|C\int_0^\varepsilon x^{a+s-2}\ln^{n-1} xdx|$, where

$T_n(x)=\sum_{j=1}^n C_n^j \frac{(2a-1)^j x^{j-1}}{(j-1)!}$ and $C$ is an appropriate constant dependent on $n$ and

independent on $\varepsilon$. Further,

$|\int_0^\varepsilon x^{a+s-2}\ln^{n-1} xdx|\leq \varepsilon^{a+s-1}(|\frac{1}{a+s-1}\ln^{n-1}\varepsilon|+|\frac{n-1}{(a+s-1)^2}\ln^{n-2}\varepsilon|+...|\frac{(n-1)!}{(a+s-1)^n}|)$. In

conditions of the theorem always $\mathrm{Re}(\rho+a-1)\geq \delta_0 > 0$, thus the factor $\varepsilon^{a+\rho-1}$ tends to zero at least as $\varepsilon^{\delta_0}$, that is faster than any negative power of $\ln\varepsilon$. This, together with the evident circumstance that in the conditions of the theorem the sums $\sum_\rho \frac{1}{(a+\rho-1)^n}$ are finite for all $n$, finish the proof.

## 2. Proof of the corollary

Already at this stage, comparing (1) for the case $n=1$, which is the sum $(2a-1)\sum_\rho \frac{1}{a-\rho}$, with the quite well known

$$\sum_\rho \frac{1}{a-\rho}=\frac{1}{a}+\frac{1}{a-1}+\sum_{m=1}^\infty \frac{\varsigma'}{\varsigma}(a)+\frac{1}{2}(\psi(a/2)-\ln\pi) \qquad (9),$$

valid for all values of $a$ distinct from the Riemann-function zeroes [3], we have the proof for a particular case $n=1$ of our corollary:

$$\frac{\varsigma'}{\varsigma}(a) = -\lim_{N\to\infty}(\sum_{m\leq N}\frac{\Lambda(m)}{m^a} - \frac{N^{1-a}}{1-a}), \text{ eq. (3). The estimation of the difference}$$

between $\sum_{\rho}\hat{g}_{n,a,\varepsilon}(\rho)$ and $\sum_{\rho}\hat{g}_{n,a}(\rho)$ given above proves also the relation

$$\frac{\varsigma'}{\varsigma}(a) = -\sum_{m\leq N}\frac{\Lambda(m)}{m^a} + \frac{N^{1-a}}{1-a} + O(N^{1/2+\Delta-a}).$$

To consider the general case, we will calculate

$$\sum_{\rho}(1-\left(\frac{\rho-a}{\rho+a-1}\right)^n) = \sum_{\rho}(1-\left(\frac{\rho+a-1}{\rho-a}\right)^n)$$ differently than above, viz. using the

functional equation $\xi(z) = \frac{1}{2}z(z-1)\pi^{-z/2}\Gamma(z/2)\varsigma(z)$ [3] and the relation

$$\frac{n(1-2a)}{(n-1)!}\frac{d^n}{dz^n}((z-a)^{n-1}\ln(\xi(z)))|_{z=1-a} = \sum_{\rho}(1-(\frac{\rho-a}{\rho+a-1})^n)$$ [1], which we rewrite

using the transformation $1-a \to a$ as

$$\frac{n(2a-1)}{(n-1)!}\frac{d^n}{dz^n}((z+a-1)^{n-1}\ln(\xi(z)))|_{z=a} = \sum_{\rho}(1-(\frac{\rho+a-1}{\rho-a})^n) \qquad (10).$$

Clearly, we have

$$\frac{d^n}{dz^n}((z+a-1)^{n-1}\ln(\xi(z)))|_{z=a} = \frac{d^n}{dz^n}((z+a-1)^{n-1}\ln(1/2))|_{z=a}$$
$$+\frac{d^n}{dz^n}((z+a-1)^{n-1}\ln(z-1))|_{z=a} + \frac{d^n}{dz^n}((z+a-1)^{n-1}\ln\pi^{-z/2})|_{z=a} \qquad (11),$$
$$+\frac{d^n}{dz^n}((z+a-1)^{n-1}\ln(z\Gamma(z/2)))|_{z=a} + \frac{d^n}{dz^n}((z+a-1)^{n-1}\ln(\varsigma(z)))|_{z=a}$$

and now we analyze all terms here one by one. Trivially

$\frac{d^n}{dz^n}((z+a-1)^{n-1}\ln(1/2))|_{z=a} = 0$ and $\frac{d^n}{dz^n}((z+a-1)^{n-1}\ln\pi^{-z/2})|_{z=a} = -\frac{\ln\pi}{2}n!$. We also

have by Leibnitz rule

$$\frac{d^n}{dz^n}((z+a-1)^{n-1}\ln\varsigma(z))|_{z=a} = \sum_{l=1}^{n}C_n^l\frac{(n-1)!}{(l-1)!}(2a-1)^{l-1}\frac{d^l}{dz^l}\ln(\varsigma(z))|_{z=a} \qquad (12).$$

To calculate other terms we apply the generalized Littlewood theorem in a manner similar to that used in our Refs. [1], [12, 13].

Again, for completeness we present this theorem here.

**Theorem 3 (The Generalized Littlewood theorem).** *Let C denotes the rectangle bounded by the lines $x = X_1$, $x = X_2$, $y = Y_1$, $y = Y_2$ where $X_1 < X_2$, $Y_1 < Y_2$ and let f(z) be analytic and non-zero on C and meromorphic inside it, let also g(z) be analytic on C and meromorphic inside it. Let F(z)=ln(f(z)) be the logarithm defined as follows: we start with a particular determination on $x = X_2$, and obtain the value at other points by continuous variation along y=const from $\ln(X_2 + iy)$. If, however, this path would cross a zero or pole of f(z), we take F(z) to be $F(z \pm i0)$ according as we approach the path from above or below. Let also $\tilde{F}(z) = \ln(f(z))$ be the logarithm defined by continuous variation along any smooth curve fully lying inside the contour which avoids all poles and zeroes of f(z) and starts from the same particular determination on $x = X_2$. Suppose also that the poles and zeroes of the functions f(z), g(z) do not coincide.*

*Then*

$$\int_C F(z)g(z)dz = 2\pi i \left( \sum_{\rho_g} res(g(\rho_g) \cdot \tilde{F}(\rho_g)) - \sum_{\rho_f^0} \int_{X_1+iY_\rho^0}^{X_\rho^0+iY_\rho^0} g(z)dz + \sum_{\rho_f^{pol}} \int_{X_1+iY_\rho^{pol}}^{X_\rho^{pol}+iY_\rho^{pol}} g(z)dz \right)$$ *where the sum is over all $\rho_g$ which are poles of the function g(z) lying inside C, all $\rho_f^0 = X_\rho^0 + iY_\rho^0$ which are zeroes of the function f(z) counted taking into account their multiplicities (that is the corresponding term is multiplied by m for a zero of the order m) and which lie inside C, and all $\rho_f^{pol} = X_\rho^{pol} + iY_\rho^{pol}$ which are poles of the function f(z) counted taking into account their*

*multiplicities and which lie inside C. The assumption is that all relevant integrals in the right hand side of the equality exist.*

To apply this Theorem, for integer $n$ greater than 1 we introduce the function $g(z) = \dfrac{n(2a-1)(z+a-1)^{n-1}}{(z-a)^{n+1}} - \dfrac{n(2a-1)}{(z-a)^2}$ and rectangular contour $C$ with vertices at $\pm X \pm iX$ with real $X \to +\infty$, and consider a contour integral $\int_C g(z)\ln(z\Gamma(z/2))dz$. (If some pole of the gamma – function occurs on the contour, just shift it a bit to avoid this). Known asymptotic of the logarithm of the gamma-function for large $|z|$, $\ln\Gamma(z) \cong O(z\ln z)$ guaranties the "disappearance" of the contour integral value (it tends to zero when $X \to \infty$ due to the asymptotic $g(z) \cong O(1/z^3)$) thus we get, after division by $2\pi i$ and elementary transformations

$$\dfrac{1}{(n-1)!}\dfrac{d^n}{dz^n}((z+a-1)^{n-1}\ln(z\Gamma(z/2)))|_{z=a} =$$
$$\dfrac{1}{n(2a-1)}\sum_{k=1}^{\infty}(1-\left(1+\dfrac{2a-1}{-2k-a}\right)^n - \dfrac{n(2a-1)}{-2k-a}) + \dfrac{1}{2}\psi\left(\dfrac{a}{2}\right) + \dfrac{1}{a} \quad (13).$$

Here $\psi(z)$ is a digamma function. In the sum occurring in this expression one easily recognizes the sum $\sum_{\rho_f^{pol}}^{X_{\rho}^{pol}+iY_{\rho}^{pol}} \int_{X_1+iY_{\rho}^{pol}} g(z)dz$ taken over simple poles $-2k$, $k=1, 2, 3\ldots$, of the gamma function (pole at $z=0$ is cancelled by the factor $z$). Clearly, $\dfrac{n(2a-1)(z+a-1)^{n-1}}{(z-b)^{n+1}} = \dfrac{d}{dz}(1-\left(\dfrac{z+a-1}{z-a}\right)^n)$ which explains why function $g(z)$ is used here: the term $-\dfrac{n(2a-1)}{(z-a)^2}$ is added just to ensure the asymptotic $g(z) \cong O(1/z^3)$ necessary to bring the contour integral value to zero. The summands $\frac{1}{2}\psi\left(\frac{a}{2}\right)+\frac{1}{a}$ in eq. (13) come exactly from this additional term.

For the summ occuring in (13), we have

$$-\sum_{k=1}^{\infty}(1-\left(1+\frac{2a-1}{-2k-a}\right)^{n}-\frac{n(2a-1)}{-2k-a}) = -\sum_{k=1}^{\infty}(1-\left(1+\frac{1-2a}{2k+a}\right)^{n}+\frac{n(2a-1)}{2k+a})$$

$$=\sum_{k=1}^{\infty}\sum_{l=2}^{n}C_{n}^{l}\frac{(1-2a)^{l}}{(2k+a)^{l}} = \sum_{l=2}^{n}C_{n}^{l}2^{-l}(1-2a)^{l}\sum_{k=1}^{\infty}\frac{1}{(k+a/2)^{l}}$$

$$=\sum_{l=2}^{n}C_{n}^{l}2^{-l}(1-2a)^{l}(\varsigma(l,\ a/2)-\frac{2^{l}}{a^{l}}) = \sum_{l=2}^{n}C_{n}^{l}2^{-l}(1-2a)^{l}\varsigma(l,\ a/2)+1-(-1+\frac{1}{a})^{l}.$$ The

change of the summation order is certainly justified here.

Similarly, using the function $\tilde{g}(z) = \frac{n(2a-1)(z+a-1)^{n-1}}{(z-a)^{n+1}}$ and the same contour, generalized Littlewood theorem gives

$$\frac{1}{(n-1)!}\frac{d^{n}}{dz^{n}}((z+a-1)^{n-1}\ln(z-1))|_{z=a} = \frac{1}{n(2a-1)}(1-(1+\frac{1}{1-a})^{n}) \qquad (14).$$

Collecting everything together and using (12), we have

$$\sum_{\rho}(1-\left(1-\frac{\rho-a}{\rho+a-1}\right)^{n}) = \sum_{\rho}(1-\left(1-\frac{\rho+a-1}{\rho-a}\right)^{n}) = 2-(-1+\frac{1}{a})^{n}-(-1+\frac{1}{1-a})^{n}+$$

$$\sum_{j=1}^{n}C_{n}^{j}(2a-1)^{j}\left\{\frac{1}{(j-1)!}\frac{d^{j}}{dz^{j}}\ln\varsigma(z)|_{z=a}\right\}+\frac{n}{2}(2a-1)(\psi(a/2)-\ln\pi)+ \qquad (15),$$

$$\sum_{j=2}^{n}C_{n}^{j}(-1)^{j}2^{-j}(2a-1)^{j}\varsigma(j,a/2)$$

which is unconditionally true for any $a$ distinct from the Riemann function zeros. Comparison of eqs. (1) and (15) finishes the proof of the Corollary; the estimations of the difference between $\sum_{\rho}\hat{g}_{n,a,\varepsilon}(\rho)$ and $\sum_{\rho}\hat{g}_{n,a}(\rho)$ given above prove the $O$-terms occurring there.

**Remark 3.** Without any change we can prove the following minor *unconditional* theorem:

**Theorem 4**. *For an arbitrary complex $a=1+it, t\neq 0$, we have:*

$$\sum_{\rho}(1-\left(\frac{\rho-a}{\rho+a-1}\right)^n) = \sum_{\rho}(1-\left(\frac{\rho+a-1}{\rho-a}\right)^n) = 2-(-1+\frac{1}{a})^n-(-1+\frac{1}{1-a})^n$$

$$\sum_{j=1}^{n} C_n^j (2a-1)^j \left\{ \frac{(-1)^j}{(j-1)!} \lim_{N\to\infty} (\sum_{m\le N} \frac{\Lambda(m)\ln^{j-1} m}{m^a} - \int_0^N x^{-a}\ln^{j-1} x\, dx) \right\} +$$

$$+\frac{n}{2}(2a-1)(\psi(a/2)-\ln\pi) + \sum_{j=2}^{n} C_n^j (-1)^j 2^{-j}(2a-1)^j \varsigma(j, a/2)$$

The fact that now $\lim_{\varepsilon\to 0+}\sum_{\rho} \hat{g}_{n,a,\varepsilon}(\rho) = \sum_{\rho} \hat{g}_{n,a}(\rho)$ is proven by *mutatis mutandi* repetition of what is said by Bombieri and Lagarias during the proof of their Theorem 2 [4]. Indeed, formula $\sum_{m\le N}\frac{\Lambda(m)}{m^{1+it}} + \frac{\varsigma'}{\varsigma}(1+it) - i\frac{N^{-it}}{t} = o(1)$ is obtained in Titchmarsh book [3], see paragraph 3.14, by other method.

The case $a=1$ also can be analyzed along similar lines. Now this is of course impossible to consider $\frac{d^n}{ds^n}\ln\varsigma(s)|_{s=1}$ but $\frac{d^n}{ds^n}(\ln(\varsigma(s)\cdot(s-1))|_{s=1}$ should be considered instead. In doing so, we get

$$\frac{d^n}{ds^n}(\ln(\varsigma(s)\cdot(s-1))|_{s=1} = (-1)^n \lim_{N\to\infty}(\sum_{m\le N}\frac{\Lambda(m)\ln^{n-1} m}{m} - \frac{\ln^n N}{n}) \qquad (16).$$

Of course, this coincides with the result of Ref. [4] where it is shown that

$$\ln(s\cdot\varsigma(s+1)) = -\sum_{n=0}^{\infty}\eta_n \frac{s^{n+1}}{n+1} \text{ where } \eta_n = \frac{(-1)^n}{n!}\lim_{N\to\infty}(\sum_{m\le N}\frac{\Lambda(m)\ln^n m}{m} - \frac{\ln^{n+1} N}{n+1}).$$

### 3. Numerical results.

Here we present some numerical calculations to illustrate the formulae exposed in the previous paragraphs in particular case of $n=1$, i.e. equation (3). All the calculations and the graphical representations were performed using the computer algebra system *Mathematica*.

Figure 1 presents the graph of the function $\frac{\zeta'}{\zeta}(a)$ compared to the term $-\sum_{m=2}^{N}\frac{\Lambda(m)}{m^a}+\frac{N^{1-a}}{1-a}$ for $\frac{1}{2}<a\leq\frac{3}{4}$. The calculation was carried out setting *N=1'000'000* in equation (3).

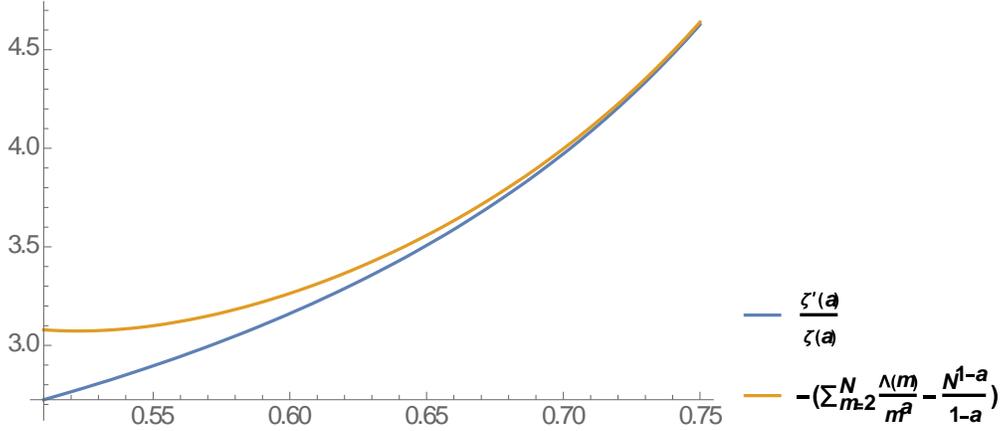

**Figure 1**: $\frac{\zeta'}{\zeta}(a)$ vs. $-\sum_{m\leq N}\frac{\Lambda(m)}{m^a}+\frac{N^{1-a}}{1-a}$

In Figure 2 we compare the trend of the function $\left|\frac{\zeta'}{\zeta}(a)+\sum_{m\leq N}\frac{\Lambda(m)}{m^a}-\frac{N^{1-a}}{1-a}\right|$ to the "control" term $O(N^{1/2-a})$.

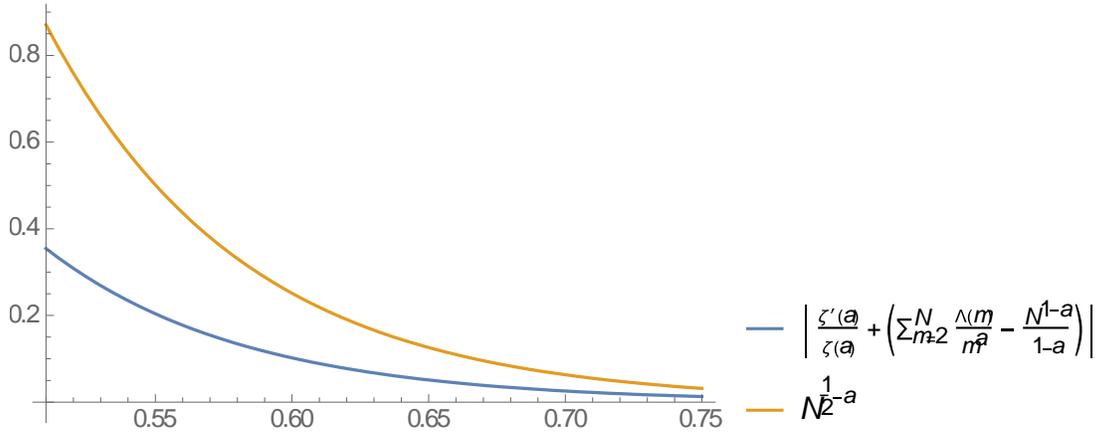

**Figure 2**: $\left|\frac{\zeta'}{\zeta}(a)+\sum_{m\leq N}\frac{\Lambda(m)}{m^a}-\frac{N^{1-a}}{1-a}\right|$ vs. $N^{\frac{1}{2}-a}$

In the last two figures, Figure 3 and 4, we analyze the convergence until $N=1'000'000$ of equation (3) for two points: one ($a = 0.55$) near the critical line and the other one ($a = 0.95$) near the right border of the critical strip.

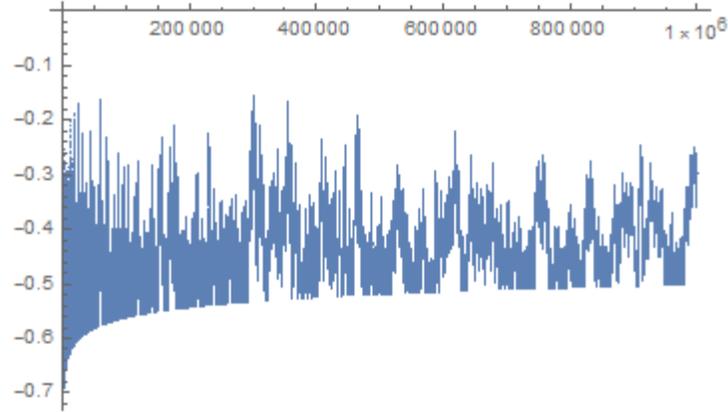

**Figure 3.** Convergence of $\left| \dfrac{\varsigma'}{\varsigma}(a) + \sum_{m \leq N} \dfrac{\Lambda(m)}{m^a} - \dfrac{N^{1-a}}{1-a} \right| - N^{1/2-a}$ at $a = 0.55$.

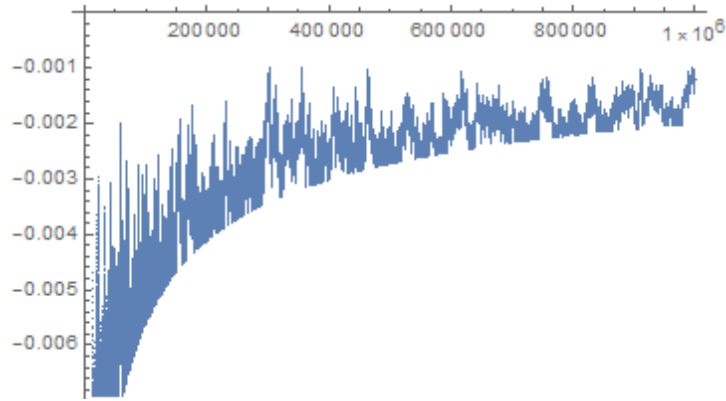

**Figure 4.** Convergence of $\left| \dfrac{\varsigma'}{\varsigma}(a) + \sum_{m \leq N} \dfrac{\Lambda(m)}{m^a} - \dfrac{N^{1-a}}{1-a} \right| - N^{1/2-a}$ at $a = 0.95$

From these data we see that while the convergence for $a=0.95$ is reasonably fast, that near the critical line is very slow (actually we are even not in a position to claim convergence here at all). This is not surprising because the argument under the *O*-sign close to the critical line is quite large, for our case $N^{-0.05} = 10^{-0.3} \cong 0.5$. (It is instructive to compare this with $N^{-0.45} = 10^{-2.7} \cong 0.002$ pertinent for *a=0.95*).

## 4. Conclusion

Thus, assuming the non-vanishing of the Riemann zeta-function in certain sub-strip of the critical strip, in the present paper we have established the formulae to calculate the value of the derivatives $\frac{d^n}{ds^n}\ln\varsigma(s)$ there, see Theorem 1 and its corollary. Interesting, and to certain extent probably even surprising, is the circumstance that the "compensating term" $\frac{N^{1-a}}{1-a}$, occurring here for $\varsigma'/\varsigma$, is exactly the same as such term occurring for the (unconditional) calculation of the Riemann function itself: $\varsigma(a) = \sum_{m \leq N}\frac{1}{m^a} - \frac{N^{1-a}}{1-a} + O(N^{-\operatorname{Re} a})$ [3]. Certainly, qualitatively this can be understood, for $\sum_{m=1}^{N}\frac{1}{m^a} \cong \int_1^N x^{-a}dx$ while $\sum_{m=1}^{N}\frac{\Lambda(m)}{m^a} \cong \int_1^N \frac{\ln x}{\ln x}x^{-a}dx \cong \int_1^N x^{-a}dx$ also - here logarithm in the nominator comes from the "magnitude" of Mangoldt function and that in the denominator from the "density of primes". Still such an exact compensation generally speaking should not occur and, in our opinion, this reflects some deep and not yet understood properties of the Riemann $\varsigma$–function.

We hope that the obtained results might find interesting and important applications in the number theory.